\documentclass{amsart}
\usepackage{amssymb}
\usepackage{amsmath}
\usepackage{amsthm}
\usepackage{tikz}
\usepackage{xcolor}
\usepackage{colortbl}
\usepackage[T1]{fontenc}
\usepackage[backref=page]{hyperref}
\hypersetup{
    colorlinks,
    linkcolor={red!80!black},
    citecolor={violet},
    urlcolor={blue!80!black}
}
\usetikzlibrary{matrix,arrows, cd, calc}  
\theoremstyle{definition}

\newtheorem{problem}{Problem}
\newtheorem{closedproblem}[problem]{(Solved) problem}

\DeclareMathOperator{\Spec}{Spec}

\DeclareMathOperator{\Hilb}{Hilb}
\DeclareMathOperator{\Slip}{Slip}
\DeclareMathOperator{\Quot}{Quot}
\DeclareMathOperator{\tr}{trace}%
\DeclareMathOperator{\Sym}{Sym}%
\DeclareMathOperator{\Ann}{Ann}%
\DeclareMathOperator{\Gr}{Gr}%
\DeclareMathOperator{\gr}{gr}%

\newcommand{\onto}{\twoheadrightarrow}

\newcommand{\kk}{\Bbbk}%

\newcommand{\BBname}{Bia{\l}ynicki-Birula}
\newcommand{\Hilbnd}{\Hilb_d(\mathbb{A}^n)}%
\newcommand{\Quotrnd}{\Quot_d^r(\mathbb{A}^n)}%
\newcommand{\Quotrdthree}{\Quot_d^r(\mathbb{A}^3)}%
\newcommand{\Quotrndpt}{\Quot_d^r(\mathbb{A}^n, 0)}%
\newcommand{\Quotrndpthree}{\Quot_d^r(\mathbb{A}^3, 0)}%
\newcommand{\lci}{\operatorname{lci}}%
\newcommand{\Gor}{\operatorname{Gor}}%
\newcommand{\sm}{\operatorname{sm}}%
\newcommand{\OO}{\mathcal{O}}%
\newcommand{\mm}{\mathfrak{m}}%

\begin{document}

\title[Open problems in deformations of Artinian\ldots]{Open problems in deformations of Artinian algebras, Hilbert schemes and
around}
\author{Joachim Jelisiejew}
\address{Institute of Mathematics, University of Warsaw, Banacha 2, 02-097
Warsaw, Poland} 
\email{j.jelisiejew@uw.edu.pl}
\thanks{Partially supported by National Science Centre grant
2020/39/D/ST1/00132. This work is partially supported by  the Thematic
Research Programme ``Tensors: geometry, complexity and quantum entanglement'',
University of Warsaw, Excellence Initiative – Research University and the
Simons Foundation Award No. 663281 granted to the Institute of Mathematics of
the Polish Academy of Sciences for the years 2021-2023.}

\begin{abstract} 
    We review the open problems in the theory of deformations of zero-dimensional objects,
    such as algebras, modules or tensors. We list both the well-known ones
    and some new ones that emerge from applications.
    In view of many advances
    in recent years, we can hope that all of them are in the range of current
    methods.

    Added 2026: this is a slightly updated version of the problem list. The aim is to mark the
    problems which were solved. The reference list was updated only to achieve this aim.
    Hence, many recent interesting works on Hilbert schemes are not mentioned
    here. We apologise for this.
\end{abstract}

\date{\today{}}
\subjclass[2020]{Primary 13D10, 14C05; Secondary 13H10}
\keywords{Hilbert scheme of points, Quot scheme of points, zero-dimensional
schemes, Artin algebras, deformations, tensors, algebraic geometry of complexity theory}

\maketitle

\begin{quotation}
    On sait qu'on ne sait rien sur cette sorte
    d'alg\`ebre.\flushright{Andr\'e Weil~\cite[p.221]{Weil_Collected}}.
    Quoted thanks to~\cite{Shafarevich__Degeneration}.\footnote{Even more engaged is
    Weil's polemic~\cite[footnote, p.221]{Weil_Collected} with Weyl~\cite[p.711]{Weyl_reactionary}.}
\end{quotation}

\section{Introduction}

    The Hilbert scheme of points of a variety $X$
    (see, e.g.,~\cite{Bertin__punctual_Hilbert_schemes, HarDeform,
    Sernesi__Deformations, Stromme_Hilbert}), occupies a special place among the moduli
    spaces investigated in algebraic geometry. On the one hand, for
    smooth curves or surfaces $X$ it behaves very nicely and gives rise to
    beautiful global geometry~\cite{Beauville__Hilbert_K3, Haiman_macdonald,
    Haiman_vanishing2}.
    On the other hand, for a smooth variety
    $X$ with $\dim X \geq 3$ it presents a mixture of pathological and unknown behaviour which
    would almost surely doom it to oblivion if not for the
    applications.

    \subsection{Applications}
    Since the Hilbert scheme parameterizes relatively elementary objects: tuples
    of points, finite-dimensional algebras, etc., deep and important applications arise
    naturally even in mathematical areas far from algebra, such as motivic
    homotopy
    theory~\cite[\S3.5]{elmanto_hoyois_Khan_Sosilo_Yakerson__Module_over_algebraic_cobordism},
    string theory~\cite{Gorsky_Nekrasov_Rubtsov} or classical
    topology~\cite{Michalek}.

    From the personal perspective of the author, a less known but key source
    of applications, probably for many years to come, is the problem of finding
    new bounds on the asymptotic complexity of matrix multiplication.
    Finding this complexity is one of the most important problems in the whole
    field of computational complexity~\cite[\S15]{BurgisserBook}, only a
    little less important than the famous ``P = NP'' problem.
    Roughly
    speaking, the complexity is
    $n^{\omega}$, where $\omega$ is an unknown constant $2\leq \omega<
    2.373$~\cite{CW, Virgi, LeGall}.
    The leading algorithm for computing $\omega$ is called the laser
    method~\cite[\S15.8]{BurgisserBook}. The
    input of this algorithm is a minimal border rank tensor. The main source
    of such tensors are multiplication tensors of \emph{smoothable} algebras
    and of modules in the principal component of
    Quot~\cite{Blaser_Lysikov, Landsberg_Michalek__Abelian_Tensors,
    Jelisiejew_Landsberg_Pal__Concise_tensors_of_minimal_brank}. The most
    successful inputs
    so far are the \emph{big Coppersmith-Winograd} tensors, which are indexed
    by $q$; to every $q$ we associate
    the multiplication tensor of the apolar algebra of $x_1^2 +
    \ldots +
    x_{q}^2$. See~\cite{Homs_Jelisiejew_Michalek_Seynnaeve} for examples of
    bounds from other tensors.
    Understanding smoothable algebras is thus understanding the possible
    inputs for the laser method.

    \newcommand{\Abarn}{\overline{A}_{n}}%
    \newcommand{\sMn}{sM_{\langle n\rangle}}

    \subsection{Big picture} The geometry of the Hilbert scheme is so rich
    that it is easy to lose sight of the big picture. We do not plan to
    discuss it here too much. However, we mention
    two fundamental invariants. The
    first one is the degree: the intuition is
    that analysing the objects of degree $d$ is roughly an exponential amount
    of work with respect to $d$. Consequently, the field divides into the
    ``near cosmos'' of small enough $d$ and the vastness of general $d$. The
    second comes from the dimension $n$ of the ambient space $X$, which we
    assume to be $\mathbb{A}^n$. This subdivides the
    spaces into the ``nicely behaved'', such as $\Hilb_d(\mathbb{A}^2)$, the
    ``border'' ones, such as $\Hilb_d(\mathbb{A}^3)$, and the ``bad'' ones,
    such as $\Hilb_d(\mathbb{A}^n)$, $n\geq 4$, see Table~\ref{tab:components}
    below. Apart from the Hilbert scheme, there is the Quot scheme, with
    similar behavior, and other moduli spaces such as the nested
    Hilbert schemes. There is also the natural Bilinear scheme, which we
    introduce in~\S\ref{sec:bilinear}, of which
    Hilbert and Quot scheme are special cases. Finally, the deformation theory
    is tightly tied to the theory of tensors. Regretfully, a review of connections of tensors and
    deformation theory is yet to come, see~\cite{Landsberg__tensors,
    Jelisiejew_Landsberg_Pal__Concise_tensors_of_minimal_brank} for some
    ideas on how it looks like. Below we keep the tensor
    language to a minimum, see~\S\ref{sec:bilinear} for some discussion.

    \subsection{Definitions and notation}\label{sec:definition}

    This section contains the bare minimum of definitions necessary to avoid
    ambivalence in the formulations of the problems below. For readers
    willing to learn more about the moduli spaces, we recommend the
    sources~\cite{Bertin__punctual_Hilbert_schemes, HarDeform,
    Sernesi__Deformations, Stromme_Hilbert} mentioned already above.

    We work over an algebraically closed field $\kk$.
    We say that an algebra $A$ or a module $M$ is \emph{finite of degree $d$}
    if $\dim_{\kk} A = d$ or $\dim_{\kk} M = d$, respectively. It is also
    common to use the name \emph{length} or \emph{multiplicity} for what we call degree.
    A finite algebra $A$ (respectively, a module $M$) is a product of local
    algebras $A_{\mm}$ (respectively, a direct sum of localizations
    $M_{\mm}$).

    Take a local algebra $A$ with the maximal ideal $\mm$. The map $\kk\to
    A/\mm$ is an isomorphism since $\kk$ is algebraically closed. The
    \emph{(local) Hilbert function} of $A$ is given by $H_A(i) :=
    \dim_{\kk} \mm^i/\mm^{i+1}$.  We have $H_A(i) = 0$ for $i$ large enough and so we
    frequently write the Hilbert function of $A$ as a vector of its nonzero
    values, for example for $A = \kk[x, y]/(x^2, x - y^2)$ we have $H_A(0) =
    H_A(1) = H_A(2) = H_A(3) = 1$, $H_A(4) =0$, so the Hilbert function is
    $(1,1,1,1)$. The \emph{associated graded algebra of $A$} is
    $\gr A := \bigoplus_{i\geq 0} \mm^i/\mm^{i+1}$.
    The Hilbert function $H_A$ is arguably the most important numerical invariant of $A$.
    Below (for example after Problem~\ref{prob:smallDegreeComponents}) we
    sometimes informally refer to components and loci by describing the Hilbert function of a
    general point, for example the ``$(1,7,7)$ component'' is a shorthand for
    ``the (unique/known/referred to) component whose general point corresponds to a local
    algebra with Hilbert function $(1,7,7)$''.

    For a
    quasi-projective scheme $X$ over $\kk$, let
    $\Hilb_d(X)$ denote the Hilbert scheme of $d$ points on $X$, for
    constructions see~\cite{fantechi_et_al_fundamental_ag,
    Haiman_Sturmfels__multigraded, Bertin__punctual_Hilbert_schemes}. More
    generally, let $\Quot_d^{r}(X)$ denote the Quot scheme whose points are,
    up to equivalence, surjections $\OO_X^{\oplus r}\onto \mathcal{M}$ of
    $\OO_X$-modules, where $\mathcal{M}$ is finite
    and has degree $d$. We will be mostly
    interested in the case when $X = \mathbb{A}^n$ and here the points of
    $\Quotrnd$ are, up to equivalence, surjections $S^{\oplus r}\onto M$ of $S$-modules, where $S =
    \kk[x_1, \ldots ,x_n]$ and $M$ is finite of degree $d$.
    The ``up to equivalence'' above just means that a surjection $\pi\colon S^{\oplus r}\onto M$ is
    identified with $S^{\oplus r}\onto S^{\oplus r}/\ker \pi$.
    The Hilbert scheme is equal to the Quot scheme taken for $r=1$ and frequently it
    is useful to think about it this way.

    The \emph{Gorenstein locus} $\Hilb_d^{\Gor}(\mathbb{A}^n)$ is an
    open~\cite[Proposition~2.5]{Casnati_Notari_6points}
    locus parameterizing Gorenstein algebras.
    The Hilbert-Chow
    map~\cite[Theorem~2.16]{Bertin__punctual_Hilbert_schemes} is a
    projective morphism
    \[
        \rho_{X,d}\colon \Hilb_d(X)\to \Sym^d(X)
    \]
    that maps a subscheme $\Gamma$ to its
    support counted with multiplicities.  The \emph{punctual Hilbert scheme}\footnote{Caution: some authors use the term
    ``punctual Hilbert scheme'' for what we call here ``the Hilbert scheme of
points''.}
    $\Hilb_d(X, x)$ is by definition $\rho_{X,d}^{-1}(d[x])$. This is a
    projective scheme and its set of $\kk$-points is the set of $\Gamma$
    supported only at $x$. The Hilbert-Chow morphism naturally generalizes to the
    ``Quot-Chow'' morphism and in particular yields a \emph{punctual Quot scheme}
    $\Quotrndpt$.

    A variant of the Hilbert scheme is the \emph{nested} or \emph{flag Hilbert
    scheme} $\Hilb^{d_1, \ldots ,d_k}(X)$. Its $\kk$-points are tuples
    \[
        \Gamma_1\subseteq \Gamma_2 \subseteq  \ldots \subseteq \Gamma_k
        \subseteq X
    \]
    where each $\Gamma_i$ is finite of degree $d_i$. The nested Hilbert scheme
    is a closed subscheme of the product of $\Hilb_{d_i}(X)$,
    see~\cite[\S4.5]{Sernesi__Deformations}, \cite{Cheah__thesis}.

    The Hilbert and Quot schemes exist also for schemes over
    $\Spec(\mathbb{Z})$, in particular
    $\Hilb_d(\mathbb{A}^n_{\mathbb{Z}})$ is a quasi-projective scheme over
    $\Spec(\mathbb{Z})$ and the fiber of
    $\Hilb_d(\mathbb{A}^n_{\mathbb{Z}})\to \Spec(\mathbb{Z})$ over every prime $p\in \Spec(\mathbb{Z})$ is
    $\Hilb_d(\mathbb{A}^n_{\mathbb{Z}/p})$, while the fiber over the generic
    point $(0)\in \Spec(\mathbb{Z})$ is $\Hilb_d(\mathbb{A}^n_{\mathbb{Q}})$.

    Below we sometimes use Macaulay's inverse systems (also known as
    apolarity). They are a useful tool for constructing finite algebras and
    modules. See for
    example~\cite{iarrobino_kanev_book_Gorenstein_algebras,
    EliasRossiShortGorenstein, Jel_classifying} for an introduction.

    An irreducible component of a finite type scheme is \emph{generically nonreduced} if it is nonreduced everywhere,
    that is, all its local rings are (geometrically) nonreduced. Otherwise is
    it \emph{generically reduced}. Being generically reduced is the same as
    having a smooth point. A finite type scheme is \emph{generically
    reduced} if all its irreducible components are.

    \subsection{Acknowledgements} The sources~\cite{iarrobino_10years,
    iarrobino_kanev_book_Gorenstein_algebras, Peeva_Stillman, aimpl} were very useful
    for collecting the more classical problems. Equally important were the
    suggestions of Jaros\l{}aw Buczy{\'n}ski, {\'A}d{\'a}m Gyenge, Mateusz Micha\l{}ek, Ritvik
    Ramkumar, Maria Evelina Rossi, Alessio Sammartano, Klemen \v{S}ivic, and
    Bal\'azs Szendr{\H{o}}i for which I am very grateful. Finally, I thank the editors of
    the CONM volume: Anthony Iarrobino, Pedro Macias Marques, Maria Evelina
    Rossi and Jean Vall\`es for the possibility to contribute. I
    also thank two anonymous referees, especially one of them, for numerous
    most valuable remarks and suggestions which improved the presentation.

    \section{Open problems}

    Below, we provide a list of open problems (coincidentally, one per each year
    since the conference that yielded the
    overview~\cite{iarrobino_10years}!).
    We hope it does list most of the problems
    ``known to experts'', but it also includes a number of directions which we
    believe are very important for the future even though they do not appear much in the
    literature yet.  We begin with two tantalizing classical
    problems, important also for applications, for example to the theory of
    counting on threefolds.

    \begin{problem}\label{problem:nonreducedThreefold}
        Is $\Hilb_d(\mathbb{A}^3)$ nonreduced for some $d$?
    \end{problem}
    \begin{problem}\label{problem:QuotSurface}
        Is the Quot scheme $\Quot_d^r(\mathbb{A}^2)$ reduced
        for all $r$, $d$?
    \end{problem}
    Note that the Quot scheme for a nontrivial sheaf can be
    nonreduced~\cite{Kass}. The preprint~\cite{Charbonnel} claims to prove
    that the answer to Problem~\ref{problem:QuotSurface} is affirmative. As far as the author knows, it is not yet
    generally accepted but neither any particular issue is identified.
    A more recent preprint~\cite{Sheshmani} also implies that the answer is
    affirmative, however there is a gap in the argument of the current arXiv
    version (from May 2025), as the authors acknowledge.
    The
    preprint~\cite{Stark} conjectures that not only $\Quot_d^r(\mathbb{A}^2)$
    is reduced, but it has rational singularities. For $r=1$ the
    answer is affirmative, as $\Hilb_d(\mathbb{A}^2)$ is
    smooth~\cite{fogarty}.

    \subsection{Irreducible components}\label{sec:component}

    Each of the moduli spaces defined above has a distinguished ``nice''
    open locus whose closure is the ``main component'':
    \begin{itemize}
        \item for the Hilbert scheme $\Hilb_d(X)$ this is the locus of $\Gamma
            \subseteq X$ with $\Gamma$ smooth. The closure is the
            \emph{smoothable component} $\Hilb_d^{\sm}(X)$.
            A smooth
            $\kk$-scheme $\Gamma$ is isomorphic to $d$ copies of $\Spec(\kk)$, so it is a
            tuple of $d$ points of $X$.
            The dimension of this component is $d\cdot \dim X$.
        \item for the nested Hilbert scheme $\Hilb^{d_1, \ldots ,d_k}(X)$ this
            is the locus of $\Gamma_{1} \subseteq \ldots \subseteq
            \Gamma_{d_k}\subseteq X$ with every $\Gamma_i$ smooth (it is
            enough to assume that the largest one is smooth). Again, this locus
            parameterizes tuples of $d_k$ points of $X$ with additional data
            coming from selecting $d_{k-1}$, $d_{k-2}$ etc. points. The
            dimension of this locus is $d_k\cdot \dim X$.
        \item for the Quot scheme $\Quotrnd$ this is the locus of surjections $S^{\oplus r}\onto
            M$ with $M \simeq \bigoplus_{i=1}^d S/\mm_i$ with $\mm_i$ maximal;
            these are the \emph{semisimple} $S$-modules. The closure is called
            the \emph{principal component} of $\Quotrnd$. Its dimension is
            $d(n+r-1)$, see for example~\cite[\S3.4]{jelisiejew_sivic}.

            Observe that an algebra
            $S/I$ corresponds to a semisimple $S$-module exactly when
            $\Spec(S/I)$ is a tuple of $d$ points of $\mathbb{A}^n$, so the above loci
            for Hilb and Quot schemes yield the same type of component.
        \item for the punctual Hilbert scheme $\Hilb_d(X, x)$ this is the
            locus of $\Gamma \subseteq X$ with $\Gamma \simeq \Spec
            \kk[\varepsilon]/(\varepsilon^d)$ as abstract schemes, called the
            \emph{curvilinear locus}. The closure is sometimes called the \emph{curvilinear
            component}.
            The dimension is $(d-1)(\dim X - 1)$ for $X$ smooth.
    \end{itemize}
    In an ideal world, we would like to look at a main component only,
    however in the cases above we have little control over this component: we
    usually do not even understand the tangent space to the component at a
    given point! Therefore, the current gold standard is to understand when
    the main component is the whole space and, when this is not the case, find
    and describe the other components.
    In the case of Hilbert/Quot schemes, we do not need to directly consider all other
    components. By abstract nonsense, each component is obtained from
    elementary
    components~\cite[Lemma~1]{Iarrobino__number_of_generic_singularities},
    \cite[\S4.4]{jelisiejew_sivic}. A component $\mathcal{Z}$ of
    the Hilbert (respectively, Quot) scheme is
    \emph{elementary}~\cite{Iarrobino__number_of_generic_singularities} if every point of $\mathcal{Z}$
    corresponds to an algebra (respectively, a module) supported only
    at one point, see below for examples.

    The current state of the art is described in Table~\ref{tab:components} below.
    \emph{None} means that the main component is the whole space,
    \emph{exist} means that we know that the space is reducible (for some
    $d$) but we do not have a description of any component, while
    \emph{known} means that we do have a description of general points of some
    components. As of today, in all cases (except Gorenstein in codimension
    $4$) where one
    component is known, many more are known as well for perhaps larger $d$.

    If no reference is given, this means that the result follows from results
    in other rows. A more thorough discussion is given below the table.

    The table shows a curious \colorbox{pink}{``dividing wall''} of unknown but existing
    components. In the years since publication, it gave rise
    to a web of conjectures and results by Jelisiejew-Ramkumar-Sammartano.
    Unfortunately, these are not published.
    \begin{table}[h]
        \centering
        \begin{tabular}{c c c c c}
            & $n=1$ & $n=2$ & $n=3$ & $n\geq 4$\\
            Gorenstein locus & none & none
            &none~\cite{kleppe_smoothability_in_codim_three,
            kleppe__unobstructedness} &
            \cellcolor{pink}exist~\cite{emsalem_iarrobino_small_tangent_space,
            bubu2010}\\
            Hilbert scheme & none & none~\cite{fogarty} &
            \cellcolor{pink}exist~\cite{iarrobino_reducibility, iarrobino_compressed_artin} &
            known~\cite{CEVV}\\
            punctual Hilbert scheme & none & none~\cite{briancon,
            Ellingsrud_Stromme__On_the_homology} & \cellcolor{pink}exist &
            known\\
            Quot scheme & none &none~\cite{MT} & \cellcolor{pink}exist
            &known~\cite{jelisiejew_sivic}\\
            punctual Quot scheme & none &none~\cite{Baranovsky,
            Ellingsrud_Lehn} & \cellcolor{pink}exist~\cite[7.10.5]{OMeara} &known\\
            flag Hilbert scheme & none &
            \cellcolor{pink}exist~\cite{Ryan_Taylor__nested_Hilbert_schemes} &
            \cellcolor{pink}exist & known\\
        \end{tabular}
        \caption{Nontrivial components of moduli spaces. Exist/known refer to
        large enough $d$.}
        \label{tab:components}
    \end{table}

    With respect to $d$, the known reducibility statements (in characteristic
    zero) are as follows.
    By~\cite{casnati_notari_irreducibility_Gorenstein_degree_11, cjn13} the
    Gorenstein locus is irreducible for $d\leq 14$, $n\leq 5$ as well as for $d\leq 13$,
    any $n$. This locus is reducible for $d \geq 14$, $n\geq 6$
    by~\cite{emsalem_iarrobino_small_tangent_space}. It is reducible
    also for $n = 4$, $d\geq 140$, $n=5$, $d\geq 42$~\cite[\S6]{bubu2010} by a
    parameter count using compressed algebras~\cite{iarrobino_compressed_artin}.
    Actually, for $n=4$ a single component is known, with $d=140$; see comment
    below Problem~\ref{prob:oddDegree}.

    For $n\geq 4$, the Hilbert scheme is reducible precisely when $d\geq 8$,
    see~\cite{CEVV}. For $n=3$ it is irreducible for $d\leq 11$,
    see~\cite{DJNT} and reducible for $d\geq 78$,
    see~\cite[Example~4.3]{iarrobino_compressed_artin}, see
    also~\S\ref{sec:compressed}.
    The punctual Hilbert scheme for $n\geq 4$ is reducible when $d\geq
    8$ and irreducible for $d \leq 7$, see~\cite{Mazzola} and
    also~\cite[Proposition~2.7 and proof of
    Theorem~1.2]{jelisiejew_keneshlou}. For $n=3$ it is irreducible for
    $d\leq 11$~\cite[Theorem~1.3]{jelisiejew_keneshlou}, while reducible for
    $d= 18$, as the locus of compressed algebras with Hilbert function
    $(1,3,6,6,2)$ has dimension $34$, by direct check or
    \cite[Theorem~II.C]{iarrobino_compressed_artin}, and $34 = (18-1)\cdot
    (3-1)$ is equal to the dimension of
    the curvilinear locus.

    \begin{problem}
        What is the smallest $d$ for which the punctual Hilbert scheme
        $\Hilb_d(\mathbb{A}^3, 0)$ is reducible? We have $12\leq d\leq 18$.
    \end{problem}
    \begin{problem}
        Does there exist a component of $\Hilb_d(\mathbb{A}^3)$ of dimension
        smaller than $3d$?
        Does there exist a component of $\Hilb_d(\mathbb{A}^3, 0)$ of dimension
        smaller than $2(d-1)$?
    \end{problem}
    Low-dimensional components of $\Hilb_d(\mathbb{A}^3, 0)$ are of interest also for the study
    of $\Hilb_d(\mathbb{A}^3)$, as they can only occur for $d$ such that
    $\Hilb_d(\mathbb{A}^3)$ is reducible,
    see~\cite[Theorem~3.5]{Gaffney_chaining_upperbounds_onsmoothablepunctual}
    or~\cite[Proposition~2.7]{jelisiejew_keneshlou}. It seems unknown
    whether existence of a small component of
    $\Hilb_d(\mathbb{A}^3, 0)$ implies the existence of such for
    $\Hilb_d(\mathbb{A}^3)$.

    \begin{problem}
        What is the smallest $d$ such that $\Hilb_d(\mathbb{A}^3)$ is
        reducible? We have $12\leq d\leq 78$.
    \end{problem}
    While the Hilbert scheme of $\mathbb{A}^3$ was investigated for many
    years, it seems that the Gorenstein locus in $\mathbb{A}^4$ received less
    attention. It is however important for applications to secant varieties,
    see~\cite{bubu2010, jabu_ginensky_landsberg_Eisenbuds_conjecture}.
    \begin{problem}
        What is the smallest $d$ such that $\Hilb_d^{\Gor}(\mathbb{A}^4)$ is
        reducible? We have $15\leq d\leq 140$.
    \end{problem}

    It is puzzling that despite proofs of reducibility, no components other
    than the smoothable one are known (see also~\S\ref{sec:compressed})
    for $\Hilb_d(\mathbb{A}^3)$. We put
    forward the following problem.
    \begin{problem}
        Describe (an open subset of) a component of $\Hilb_d(\mathbb{A}^3)$ other than the smoothable
        one. Here $d$ can be a large integer, even $d\gg 78$.
    \end{problem}
    Finding the components on $\mathbb{A}^3$ means working with large $d$.
    However, it would be really interesting to have a full description of what
    happens ``near home''. This is the motivation for the following
    problem. (The problem was solved for $d\leq 10$
    in~\cite{Galazka_Keneshlou_Sivic__Small_components}. The only elementary
    components are the Shafarevich's ones.)
    \begin{closedproblem}\label{prob:smallDegreeComponents}
        Classify all components of $\Hilb_d(\mathbb{A}^n)$ for $d$ reasonably
        small, for example $d= 10$ or a little higher, $n$ arbitrary.
    \end{closedproblem}
    It is enough to classify elementary components: those which
    parameterize irreducible finite subschemes.
    The known elementary components for small $d$ come from
    Shafarevich~\cite{Shafarevich_Deformations_of_1de}, which
    generalizes~\cite{emsalem_iarrobino_small_tangent_space}. The components
    parameterize algebras with Hilbert function $(1,n,r)$ for certain $r$, see
    discussion below Problem~\ref{prob:ShGap}. The
        component $(1,7,7)$ was also found independently later
        in~\cite[Corollary~4.10]{bertone_cioffi_roggero__smoothable_1771}. The other
    known elementary components of degree $d\leq 15$, according to the author's
    knowledge, are the $(1,5,3,4)$
    obtained by Huibregtse~\cite{Huibregtse_elementary}, the
    $(1,5,6,1)$ obtained by Kleiman-Kleppe~\cite{Kleiman_Kleppe},
    the $(1,6,6,1)$ obtained earlier by
    Iarrobino-Emsalem by the method
    of~\cite{emsalem_iarrobino_small_tangent_space}, and the $(1,4,6,4)$
    obtained by Satriano-Staal in~\cite[Theorem~1.3]{Satriano_Staal}. It would be
    very interesting to understand whether these are all elementary components
    for small $d$. (For example, Huibregtse
    generalized his result in~\cite{Huibregtse_more_elementary}, but $d$ is
    larger in the examples obtained there). The recent advances in \BBname{}
    decompositions, see for example~\cite{Jelisiejew__Elementary}, make such a
    classification project ambitious but not insane.

    \medskip
    We now pass to discussing the Quot scheme.
    We consider the case of $r$ large enough, say $r\geq
    d$. We still assume characteristic zero. For $n\geq 4$ the scheme
    $\Quotrnd$ is reducible for every $d\geq
    4$~\cite[p.72]{Guralnick__commuting_pairs}.  For $n=3$, the scheme
    $\Quotrdthree$ is irreducible for $d\leq 10$
    \cite{Sivic__Varieties_of_commuting_matrices} and actually for $d\leq 11$
    (personal communication from Klemen \v{S}ivic), while this scheme is reducible for $d\geq
    29$, see~\cite{Holbrook_Olmadic} or~\cite[p.238]{Ngo_Sivic}. The punctual
    Quot scheme $\Quotrndpthree$ behaves
    similarly. It is irreducible for $d\leq 6$ by~\cite{Ngo_Sivic} and
    reducible for $d\geq 13$ by~\cite[7.10.5]{OMeara}.

    \begin{problem}
        What is the smallest $d$ such that $\Quotrdthree$ is
        reducible for some $r$? We have $12\leq d\leq
        29$.
    \end{problem}
    A similar question can be asked for the punctual Quot scheme.

    \goodbreak
    In the case of the flag or nested Hilbert scheme, we have an additional
    parameter: the length of the flag. Applications mostly concern the case of
    a two-element flag, see for example~\cite{Hsiao_Szenes} and references
    therein. Partially this may be because the two-element flag is used to
    construct Nakajima operators,
    see~\cite[\S8.3]{nakajima_lectures_on_Hilbert_schemes}.
    \begin{problem}
        What are the singularities of $\Hilb^{d_1,d_2}(\mathbb{A}^2)$? Can it
        be nonreduced?
    \end{problem}
    For $d_1\leq 2$ the singularities are rational,
    see~\cite{Song__Rational_sings, Ramkumar_Sammartano__Rational_sings}.
    For affine spaces with $\mathbb{A}^n$ and/or more nestings, nonreducedness can happen,
    see, for example~\cite{GGGL, GGGL_New_components, Graffeo_Lella_Components_of_nested_Hilbert_scheme_of_few_points}.

    Finally, for the Gorenstein locus we have also an additional structure:
    Iarrobino's symmetric decomposition~\cite{iarrobino_associated_graded}.
    Call a local algebra \emph{graded} if it is isomorphic to its associated
    graded $\gr A$. Iarrobino's decomposition
    yields, for a local Gorenstein algebra $A$, a canonical graded
    Gorenstein algebra $Q(0)$ with a surjection $\gr A\onto Q(0)$. This surjection
    is an isomorphism precisely when $A$ is isomorphic to its associated graded.
    The author is not aware of any examples of Gorenstein nonsmoothable
    algebras $A$ with $\gr A$ smoothable. For $H_A(4) = 0$ they do not exist
    by~\cite{CENR_Poincare}. However, they should exist in general.

    (The following problem was solved, at least in terms of showing existence, in~\cite{Wu__Qzero_existence}.)
    \begin{closedproblem}
        Find a nonsmoothable local Gorenstein algebra $A$ with $Q(0)$
        smoothable, preferably of small degree.
    \end{closedproblem}

    Many of the irreducibility results above use the characteristic zero
    assumption. However, remarkably, there are no counterexamples for these
    results in positive characteristic and, in general, the author is unaware
    of any interesting phenomena for small $d$ that exist only in positive
    characteristic.

    \subsection{Singular surfaces}

    In general, the relation between the singularities of $X$ and the
    topological properties of $\Hilb_d(X)$ for varying $d$ is complicated.
    Much is known for singular curves, we point to \cite[p.7]{Altman_Iarrobino_Kleiman}, \cite{Migliorini__HOMFLY,
    Oblomkov_Shende, Rennemo__Homology} as example references.

    Recall that by Fogarty's result, the Hilbert
    scheme of a smooth surface is smooth.
    Below we consider normal singular surfaces. By a \emph{singularity} we mean the
    complete local ring at a singular point. The nicest class of surface singularities are du Val
    singularities~\cite{Reid, Durfee__Fifteen_characterizations}. In
    algebraic terms, over an algebraically closed field $\kk$ of
    characteristic zero, a surface singularity is \emph{du Val} if
    both of the following hold
    \begin{enumerate}
        \item the tangent space is three-dimensional, so that the complete
            local ring is $\kk[[x,y,z]]/(f)$ for some local coordinates $x$,
            $y$, $z$.
        \item If we write $f = \sum_{i,j,k} a_{ijk}x^iy^jz^k$ then all of the
            conditions below hold
            \begin{itemize}
                \item there is some nonzero $a_{ijk}$ with $i+j<2$,
                \item there is some nonzero $a_{ijk}$ with $i+j+k<3$,
                \item there is some nonzero $a_{ijk}$ with $2i+j+k<4$,
                \item there is some nonzero $a_{ijk}$ with $3i+2j+k<6$.
            \end{itemize}
    \end{enumerate}
    Moreover, the characterization does not depend on the choice of local coordinates
    $x$, $y$, $z$ above.
    See~\cite[Theorem 2.1(3)]{Reid} for a sketch of proof.

    For a normal surface $S$ with all singularities du
    Val the Hilbert scheme $\Hilb_d(S)$ is irreducible
    by~\cite{Craw_Gammelgaard_Gyenge_Szendroi, Zheng__ADE}, see also \cite{Craw_intro}. In
    general, a normal surface $S$ has only finitely many singular points $p_1$, \ldots
    , $p_s$ and the
    irreducibility of $\Hilb_d(S)$ depends only on the singularities
    $\widehat{\OO}_{S,p_i}$ for $i=1, \ldots ,s$. In particular, using
    quasi-projectivity one can reduce to studying affine $S$
    with a unique singular point at the origin.

    \begin{problem}\label{prob:reducibilityOfSurfaces}
        Is $\Hilb_d(S)$ reducible for every normal surface $S$ with a
        singularity which is not a du Val singularity?
    \end{problem}
    For surfaces with a point of multiplicity $\geq 5$, the answer is affirmative, see~\cite{Graffeo_Lella_Components_of_nested_Hilbert_scheme_of_few_points}.

    \begin{problem}
        In the setting of Problem~\ref{prob:reducibilityOfSurfaces}, assume
        that the singularities of $S$ are du Val. Is $\Hilb_d(S)$ reduced for
        every $d$?
    \end{problem}
    The proof of irreducibility~\cite{Craw_Gammelgaard_Gyenge_Szendroi} does not imply
    reducedness, however it gives some
    ideas about a possible proof.
    It is enough to make the proof for the cases
    $S\subseteq \mathbb{A}^3$ listed in~\cite[Table~1, p.5]{Reid}.

        \newcommand{\VCL}{\operatorname{VC}}%

            \subsection{Very compressed algebras and modules}\label{sec:compressed}
            Compressed algebras, developed and employed successfully by
            Iarrobino~\cite{iarrobino_reducibility, iarrobino_compressed_artin,
            Kleiman_Kleppe} are the main known tool to construct large loci in the
            Hilbert scheme and in its Gorenstein locus.
            The idea is very natural and has numerous variants and
            adaptations.

            For fixed $n$ and $d$ we can form the \emph{locus of very
            compressed algebras} $\VCL_{n,d}^\circ$. It consists of degree $d$ algebras
            $A = \kk[x_1, \ldots ,x_n]/I$, where the radical of $I$ is $(x_1, \ldots
            ,x_n)$ and with Hilbert function maximal
            for ``as long as possible'' given the degree. For example
            $\VCL_{4, 21}^{\circ}$ parameterizes algebras with Hilbert function
            $(1,4,10,6)$. More formally, if we pick an $s$ so that
            $\binom{n+s-1}{n} < d \leq \binom{n+s}{n}$, then the Hilbert
            function of $A$ satisfies
            \[
                H_A(i) = \begin{cases}
                    \binom{i+n-1}{n-1} & \mbox{for } i< s,\\
                    d - \binom{n+s-1}{n} & \mbox{for } i =s,\\
                    0 &\mbox{for } i > s.
                \end{cases}
            \]
            In this case $\VCL_{n,d}^{\circ}$ is isomorphic to the Grassmannian $\Gr(d
            - \binom{n+s-1}{n}, \binom{n-1+s}{s})$.
            By moving the support, we obtain a yet bigger locus
            $\VCL_{n,d} = \VCL_{n,d}^{\circ}\times \mathbb{A}^n$.
            Already for $n= 3$ one can
            check~\cite{iarrobino_reducibility} that the dimension of
            $\VCL_{3,d}$ grows faster than $3d$ so the existence of
            such a large locus proves reducedness of $\Hilb_{d}(\mathbb{A}^3)$
            for large $d$. It turns out that
            $d = 96$ is the first value when $\dim \VCL_{3, d} \geq 3d$.

            A generalization of very compressed algebras are \emph{compressed}
            algebras~\cite{iarrobino_compressed_artin}. We do not discuss them
            at length here, mostly because a natural treatment of them would
            require introducing apolarity. We only recall one specific
            important locus constructed by Iarrobino
            in~\cite[Example~4.3]{iarrobino_compressed_artin}: it is a locus
            $\mathcal{L}_{78}\subseteq \Hilb_{78}(\mathbb{A}^3)$ which has
            dimension $\dim \mathcal{L}_{78} = 235 > 3\cdot 78$, so it proves
            that $\Hilb_{78}(\mathbb{A}^3)$ is reducible.

        The compressed algebras allow us to formulate a more specific
        incarnation of Problem~\ref{problem:nonreducedThreefold}, actually two
        slightly different ones.
        \begin{problem}\label{problem:VCLthree}
            Consider the locus $\mathcal{L}_{78}$. Is this locus a component? (If so,
            the component is generically nonreduced.)
            Consider the locus $\VCL_{3, 96}$. Is this locus a component? (If so,
            the component is generically nonreduced.)
        \end{problem}

        The
        rationale behind this is that since~\cite{iarrobino_compressed_artin}
        there has been no progress on finding large loci and it seems hard to
        imagine how to construct any locus containing
        $\mathcal{L}_{78}$ as a proper subset. For $\VCL_{3, 96}$ the
        rationale is similar, but the evidence is weaker as there may be loci of
        compressed algebras in $\Hilb_{96}(\mathbb{A}^3)$ which have dimension
        higher than $\dim \VCL_{3, 96}$.
        An even bolder version of Problem~\ref{problem:VCLthree} is the
        following.

        \begin{problem}
            Is the following alternative true for all $(n, d)$:
            \emph{$\VCL_{n, d}$ lies in the
            smoothable component or it forms a component}?
        \end{problem}
        For $n=3$, $d\leq 95$ the locus $\VCL_{n,d}$ is contained in the
        smoothable component, at least in characteristic zero~\cite{DJNT}.
        It can be checked as in~\cite[\S6.5]{jelisiejew_sivic} that for $n=4$, $d=21,22,23,24$ we
        indeed get that $\VCL_{4,d}$ is a (generically nonreduced) component of
        $\Hilb_d(\mathbb{A}^4)$; this
        solves~\cite[Problem~3.8]{aimpl}.

        A subproblem of the above is the case when $d \leq
        \binom{n+2}{2}$, so that $s \leq 2$.
        \begin{problem}[Shafarevich's gap, algebras with Hilbert function
            $(1,n,r)$]\label{prob:ShGap}
            Let $r = d - (1 + n)$ and assume that $\frac{(n-1)(n-2)}{6}+2 < r \leq
            \binom{n+1}{2}$. Decide for which $(n, d)$ the locus
            $\VCL_{n,d}$
            is a component or lies in the smoothable component.
        \end{problem}
        For $r$ outside the range above, the situation is as follows.
        Values $r=1,2$ give smoothable algebras, see~\cite[Proposition~4.10]{CEVV}.
        A choice $3 \leq r \leq \frac{(n-1)(n-2)}{6}+2$ yields a
        generically reduced elementary component
        $\VCL_{n,d}$~\cite[p.192]{Shafarevich_Deformations_of_1de}. The choice
        $r > \frac{n^2-1}{3}$ cannot lead to a
        generically reduced component~\cite[p.184]{Shafarevich_Deformations_of_1de} (but
        perhaps a nonreduced component is still possible). The case $n=5$ was solved in
        \cite{Zhao_Shafarevich}.

        One can also consider compressed and very compressed modules.
        Here, care is necessary to define the apolarity action properly,
        see~\cite[\S2]{Kunte__Gorenstein_modules_of_finite_length},
        \cite[\S4.1]{jelisiejew_sivic}, \cite[\S3]{Harkonen_Hirsch_Sturmfels}.
        We leave the precise analysis of compressed modules to future
        researchers. The definition of the very compressed locus is
        easier. Namely, we define a sublocus $(\VCL_{n,d}^{r})^{\circ}\subseteq
        \Quotrndpt$, which consists of modules $M$ with Hilbert function
            \[
                H_A(i) = \begin{cases}
                    r\cdot \binom{i+n-1}{n-1} & \mbox{for } i< s,\\
                    d - r\cdot \binom{n+s-1}{n} & \mbox{for } i =s,\\
                    0 &\mbox{for } i > s,
                \end{cases}
            \]
            where $s$ is such that $r\binom{n+s-1}{n} < d \leq
            r\binom{n+s}{n}$ and define $\VCL_{n,d}^{r}$ as
            $(\VCL_{n,d}^{r})^{\circ}\times \mathbb{A}^n \subseteq \Quotrnd$,
            where the embedding sends $([M], v)$ to $M$ translated by the
            vector $v$, so that the resulting module is supported at $v$.
        \begin{problem}
            Is the following alternative true for most $(n, d, r)$:
            \emph{$\VCL_{n, d}^r$ lies in the
            principal component or forms a component}?
        \end{problem}
        Here ``most'' means: ``all except a short finite list of exceptions''.
        By~\cite{jelisiejew_sivic}, the alternative holds for $d\leq 7$ except
        for $(n, d, r) = (4, 7, 3), (4,7,4)$; in many cases we do get a component.

    \subsection{Compressed Gorenstein algebras}

            For a fixed $n$ and $j$ we can consider the apolar algebra of
            a general form $F\in \kk[y_1, \ldots ,y_n]_{j}$ where
            $y_{\bullet}$ are variables dual to variables $x_{\bullet}$ in the
            ring $S =\kk[x_1, \ldots ,x_n]$. Taking $F$ to its apolar algebra
            yields a map $U\to \Hilb_{d}^{\Gor}(\mathbb{A}^n)$, where
            $U$ is an open subset of the affine space $\kk[y_1, \ldots
            ,y_n]_{j}$ consisting of $F$ whose apolar algebras have the largest
            degree, see~\cite{emsalem, jelisiejew_VSP, Kleiman_Kleppe}.
            See~\cite{Boij_betti_numbers} for a description of the Betti
            tables in the image.

    \begin{problem}[General forms of odd degree
            {\cite[Conjecture~6.30]{iarrobino_kanev_book_Gorenstein_algebras}}]\label{prob:oddDegree}
        Fix an odd $j$. Let $F$ be a general form of  degree $j$ in $n$
        variables, where $n\geq 6$ or $(n,j) = (5, \geq\!\! 5)$ or $(n, j) = (4,
        \geq\!\! 15)$ and moreover $(n,j)\neq (7, 3)$.
        Let $\Gamma = \Spec(S/\Ann(F))\subseteq \mathbb{A}^n$. Is it true that $\dim
        (T_{[\Gamma]})_{<0} = n$?
    \end{problem}
    If the claim is true, then $[\Gamma]$ lies on a unique elementary component
    of $\Hilb_{d}(\mathbb{A}^n)$ and is a smooth point there.
    The case $j=3$ is solved in~\cite{Szafarczyk}. It is likely that the
    method generalizes. By semicontinuity, it is enough to check that $\dim
    (T_{[\Gamma]})_{<0} = n$ for a single $F$, which can be done using a
    computer algebra system. The author verified the cases $(n,j) = (4,15),
    (5,5)$.

    \begin{problem}[Odd behaviour of even degrees]\label{prob:evenDegree}
        Fix an even $j>2$. Let $F$ be a general form of  degree $j$ in $n$
        variables.
        Let $\Gamma = \Spec(S/\Ann(F))\subseteq \mathbb{A}^n$. Is it true that
        then $[\Gamma]\in
        \Hilb_{d}(\mathbb{A}^n)$ lies on a unique elementary component
        which is everywhere nonreduced for ``most'' pairs $(j,n)$? Can one at
        least find one such pair? (Here $n$ should probably be taken large
        with respect to $j$.)
    \end{problem}
    The reason behind the different expectations in
    Problems~\ref{prob:oddDegree}-\ref{prob:evenDegree} is numerical. For $j
    = 2k$ and general $F$ let $I := \Ann(F)$ and $\Gamma := \Spec(S/I)$. Then we have $I_{\leq k} =
    0$ so that $(I^2)_{\leq 2k+1} = 0$. Moreover $I_{2k+1} = S_{2k+1}$, since
    $\deg(F) < 2k+1$. It follows that $\dim_{\kk} (I/I^2)_{2k+1} =
    \dim_{\kk} S_{2k+1}$ is large. By duality for Gorenstein algebras, the
    space $(I/I^2)_{2k+1}$ has the same dimension as the degree $-1$ part of
    the tangent space at $[\Gamma]\in \Hilb_d(\mathbb{A}^n)$,
    see for example~\cite[Proposition~2.12]{Szafarczyk}. This rules
    out the behaviour expected in the odd $j$ case. Moreover, the feeling is
    that these degree $-1$ tangent directions have no geometric ``meaning'', so for most
    cases they should not extend to higher order. If this is the case, then an argument similar
    to~\cite{Szachniewicz} proves nonreducedness.

    \subsection{Special components}

    In this section we collect the odd problems concerning components. All of
    them are classical and of interest.
    The following problem was solved in~\cite{Farkas_Pandharipande_Sammartano}. The bound on $n$
    was improved in~\cite{Wu__One_variable}.

    \begin{closedproblem}[\cite{lella_roggero_Rational_Components, aimpl}]
        Does there exist a non-rational component of $\Hilb_d(\mathbb{A}^n)$?
    \end{closedproblem}
    Recall that a component is rational if it contains an open subset
    isomorphic to an open subset of an affine space. The rationality is of
    interest also for more practical terms: some algorithms aiming at computing the
    local structure of the Hilbert scheme~\cite{Ilten, Lella}, when applied to
    a smooth point, will currently succeed in producing a universal family only
    if the component containing the point is
    rational.

    \begin{problem}
        Does there exist a component of $\Hilb_d(\mathbb{A}^n_{\mathbb{Z}})$
        that lives entirely in characteristic $p$, for small $d$?
    \end{problem}
    For $d$ a few hundred, an example is given
    in~\cite[\S5]{Jelisiejew__Pathologies}. A smaller example would be useful
    for finding varieties over characteristic $p$ nonliftable to
    characteristic zero.

    \begin{problem}
        Fix $n$.
        What is the smallest dimension $\Delta$ of a component of
        $\Hilb_d(\mathbb{A}^n, 0)$, when $d\to \infty$? Or, a bit more precisely,
        what can be said asymptotically about $1/d\cdot (nd - \Delta)$?
    \end{problem}

    There seemed to be a consensus among experts that every component of
    $\Hilb_d(\mathbb{A}^n, 0)$ has dimension at least $(n-1)(d-1)$. This was
    recently disproved by Satriano-Staal~\cite{Satriano_Staal} already for $n=4$.
    The dimension of the \emph{largest} component is known to grow
    asymptotically as $d^{2 - 2/n}$, see~\cite{Briancon_Iarrobino__dimension_of_punctual_Hilbert_scheme}. See
    also~\cite{poonen_moduli_space}.

    \begin{problem}[{\cite[p.313]{iarrobino_10years}}]\label{prob:rigid}
        Does there exist an irreducible finite $\kk$-scheme $\Gamma$ with $T^1_{\Gamma} =
        0$, other than $\Gamma
        = \Spec(\kk)$?
    \end{problem}
    Recall that such a scheme is called \emph{rigid} (to be safe, assume $\kk$
    has characteristic zero here). Here, $T^1_{\Gamma}$ is the first Schlessinger's
    functor~\cite[\S3]{HarDeform}. If $\Gamma = \Spec(S/I)$, with $S = \kk[x_1, \ldots ,x_n]$ and $d
    = \dim_{\kk} S/I$, then
    $T^1_{\Gamma}$ is the quotient of the tangent space
    $T_{[\Gamma]}\Hilbnd$ by the $S$-submodule generated by (the
    images in $T_{[\Gamma]}$ of) the $n$-dimensional space of vector
    translations, see for example~\cite[Proposition~3.10]{HarDeform}.
    It seems that for Gorenstein $\Gamma$ the answer is negative
    by~\cite[Theorem~2.3
    and Proposition~4.1]{Aleksandrov__duality}.

    \subsection{Tangent spaces} The tangent space in general does not seem to
    carry much geometric information, except that a large tangent space makes
    explicit calculations very resource-consuming. However, on the ``dividing
    wall'', the tangent space seems to carry a lot more information, at
    witnessed by the following problems.
 
    The Brian{\c{c}}on-Iarrobino conjecture below was solved in three variables in~\cite{Mackenzie_Rezaee}.
    \begin{closedproblem}[\cite{Briancon_Iarrobino__dimension_of_punctual_Hilbert_scheme}]
        Which is the maximal tangent
        space dimension to a point of $\Hilb_d(\mathbb{A}^3)$?
    \end{closedproblem}
    Brian{\c{c}}on-Iarrobino~\cite[p.544]{Briancon_Iarrobino__dimension_of_punctual_Hilbert_scheme}
    put forward a conjectural answer. The conjecture concerns every $n$, we discuss it here
    only for $n=3$. The conjecture has two versions. The first one is
    formulated only for degrees of the form $d =
    \binom{k+2}{3}$ for some natural $k$; we call it the \emph{triangular}
    version.
    This version claims that the
    maximal tangent space is attained at the point corresponding to $(x,y,z)^k$. The second,
    general conjecture was made for every $d$.
    The more general version was
    disproved and updated by Sturmfels~\cite{Sturmfels__counterexamples}.
    Sturmfels' version was subsequently disproved by
    Ramkumar-Sammartano~\cite{Ramkumar_Sammartano}. The triangular case of the original conjecture
    still stands and Ramkumar-Sammartano come close to proving it.

    The following two problems have very different formulations, but come
    from the common motivation in enumerative geometry~\cite{MNOP,
    Pandharipande}. (The following problem has been
    solved~\cite{Graffeo_Giovenzana_Lella}
    while the present paper was in review.)
    \begin{closedproblem}\label{prob:parity}
        Is it true that for every point $[\Gamma]\in \Hilb_d(\mathbb{A}^3)$ the parity
        of $d$ and of the dimension of the tangent space at $[\Gamma]$ is the same?
    \end{closedproblem}
    The answer is affirmative when $I = I(\Gamma)$ is
    monomial~\cite[Theorem~2]{MNOP} and more generally
    for ideals $I$ homogeneous with respect to the standard
    $\mathbb{N}$-grading (and even some
    nonstandard gradings), see~\cite{Ramkumar_Sammartano_parity}.
    While this list was in review, a negative answer to the question in
    Problem~\ref{prob:parity} was given by
    Giovenzana-Giovenzana-Graffeo-Lella~\cite{Graffeo_Giovenzana_Lella}; it is
    striking that their example has degree only
    $d=12$ and is a (non-standard) graded ideal.
    The story of Problem~\ref{prob:parity} is not yet finished. A weaker
    question: ``Do parities agree on an open dense subset?'' is still open
    and well motivated by the enumerative geometry,
    see~\cite[Theorem~1.2]{Ricolfi_notes}.

    The formulation of the second problem is very far from elementary, as it
    uses the Behrend function, a measure of the singularity, introduced by
    Behrend~\cite[Definition~1.4]{Behrend__function} and quite difficult to
    compute even for small examples.
    It is known that the function is equal to $(-1)^d$ for every point
    corresponding to a monomial
    ideal of degree $d$~\cite[Proposition~4.2]{Behrend__Fantechi}. See
    also~\cite[\S4]{Pandharipande_Thomas__Curve_counting} for an example of
    computation.
    (The following question was answered in the negative in~\cite{Kool_Jelisiejew_Schmierman}.)
    \begin{closedproblem}
        Is the Behrend function equal to $(-1)^d$ for every point of $\Hilb_d(\mathbb{A}^3)$?
    \end{closedproblem}
    If true, this would imply that $\Hilb_d(\mathbb{A}^3)$ is generically
    reduced, see~\cite{Ricolfi_notes}. In general, the positive answer would restrict the possible
    singularities of this Hilbert scheme. A negative answer would be
    philosophically very important for the research into counting on
    threefolds, as the multiplicities of the components would need to be taken
    into account.

    The remaining problems in this section follow from Hu's striking
    idea~\cite{Hu__singular_Hilbert_schemes} that
    the singularities of points with only a small surplus in the tangent space could
    be classified and understood. We follow Hu's decision to restrict to
    $\mathbb{A}^3$. Let $[\Gamma]$ be a singular point of the smoothable
    component of $\Hilb_d(\mathbb{A}^3)$. Then the tangent space
    $T_{[\Gamma]}\Hilb_d(\mathbb{A}^3)$ has dimension at least $3d+1$. According
    to~Problem~\ref{prob:parity} is could be at least $3d+2$.
    The smallest dimension known to be attained is $3d+6$, it occurs for
    example for $\Gamma = \Spec(\kk[x_1, x_2, x_3]/(x_1, x_2, x_3)^2)$.

    To formulate the problem below, we need the notion of smooth equivalence,
    see for example~\cite[\S2]{Jelisiejew__Pathologies}. Namely, two complete
    local rings $A$ and $B$ are \emph{smoothly equivalent} if $A[[
    x_1, \ldots ,x_{a}]]$ and $B[[x_1', \ldots ,x_{b}']]$ are isomorphic for
    some $a$ and $b$; the parentheses denote taking formal power series.
    For closed points $x\in X$, $y\in Y$ of finite type schemes $X$, $Y$ we
    say that $x\in X$, $y\in Y$ are
    \emph{smoothly equivalent} if their complete local rings are.
    This can be rephrased in geometric terms: closed points $x\in X$, $y\in Y$ are smoothly equivalent if and only if there is a
    finite type
    scheme $Z$ with a closed point $z\in Z$ and
    \emph{smooth} maps $Z\to X$, $Z\to Y$ which send $z$ to $x$,
    $y$, respectively. Hence,
    smooth equivalence is a way to compare complete local rings
    of different dimensions. Properties of a singularity, such as being
    reduced, are preserved by smooth equivalence, although this needs to be
    carefully proven for each given property,
    see for example~\cite[\href{https://stacks.math.columbia.edu/tag/033D}{Tag
    033D}]{stacks_project}.

    \begin{problem}[Mildly singular points of $\Hilb_d(\mathbb{A}^3)$,
        \cite{Hu__singular_Hilbert_schemes}]
        Let $[\Gamma]\in \Hilb_d(\mathbb{A}^3)$ be a point of the
        smoothable
        component (see~\S\ref{sec:component}) with tangent space of dimension
        $3d+6$. Is it true that $[\Gamma]\in \Hilb_d(\mathbb{A}^3)$ is smoothly
        equivalent to the cone point of a cone over the Grassmannian
        $\Gr(2, 6)$?
    \end{problem}
    For $d = 4$ this is true classically and is explained for example in~\cite{Katz__4points,
    stevens_deformations_of_singularities}.
    For $d\leq 7$ this
    is proven in~\cite[Proposition~4.30]{Hu__singular_Hilbert_schemes}.
    For classes of monomial ideals it was proven in~\cite{AJR}.

    \begin{problem}[Mildly singular points of
        $\Hilb_d^{\Gor}(\mathbb{A}^4)$]\label{prob:mildlySingCodim4}
        Let $[\Gamma]\in \Hilb_d^{\Gor}(\mathbb{A}^4)$ be a point of the smoothable
        component. Suppose that the surplus
        \[
            \dim T_{[\Gamma]} \Hilb_d(\mathbb{A}^4) - 4d
        \]
        is
        small. Can we classify singularities of such points up to smooth
        equivalence?
    \end{problem}
    For example,
    Ranestad-Schreyer~\cite[Cor~5.16]{Ranestad_Schreyer__polar_simplices}
    proved that for
    \[
        \Gamma = \Spec\left(\frac{\kk[x_1, \ldots ,x_4]}{(x_ix_j\ |\ i\neq j) +
        (x_i^2 - x_1^2\ |\ i\neq 1)}\right)
    \]
    the
    surplus is $5$, the singularity is equivalent to a cone point of a cone over the $10$-dimensional
    spinor variety. This is very surprising as the spinor variety is a very
    special one. One concrete version of Problem~\ref{prob:mildlySingCodim4} would
    be to consider difference $5$ and prove that we always obtain a spinor
    variety, up to smooth equivalence.

    \subsection{Punctual Hilbert schemes}

    The punctual Hilbert schemes (see~\S\ref{sec:definition}) naturally appear when considering the
    elementary components of the Hilbert scheme of points. Much is known about
    the case of two variables and a little for the case of three variables, see for
    example~\cite{briancon, Granger_memoir,ia_deformations_of_CI}.
    Still, it seems that many natural questions remain open.
    The punctual Hilbert scheme is important for applications to enumerative
    problems and algebraic homotopy theory. For example, tautological counts on the smoothable component
    of the Hilbert scheme reduce to the curvilinear Hilbert
    scheme~\cite{Berczi_One, Berczi_Two};
    understanding its components and upper bounds on the dimension is
    important here.

    \begin{problem}\label{prob:componentsPunct}
        Say that an irreducible component $\mathcal{Z}$ of $\Hilb_d(\mathbb{A}^n, 0)$ is
        \emph{smoothable} if $\mathcal{Z}$ is contained in the smoothable
        component of the Hilbert scheme $\Hilb_d(\mathbb{A}^n)$. What can be said about the smoothable
        components $\mathcal{Z}$? How many are there, for a given $(n, d)$?
    \end{problem}
    In contrast with the Hilbert scheme of points, as far as the author knows, there is
    no natural injective map from the components of $\Hilb_d(\mathbb{A}^n, 0)$ to
    the components of $\Hilb_{d+1}(\mathbb{A}^n, 0)$. It would be very interesting to see if the
    answer in Problem~\ref{prob:componentsPunct} can sometimes decrease
    with increasing $d$.

    (The following was solved
    in~\cite{Berczi_Svendsen} while the present paper was in review.)
    \begin{closedproblem}
        Do all monomial ideals in $\Hilb_d(\mathbb{A}^n, 0)$ lie in the
        closure of the curvilinear locus?
    \end{closedproblem}
    As a motivation, it is known classically that all monomial ideals lie in the smoothable component of
    $\Hilb_d(\mathbb{A}^n)$, see~\cite{har66connectedness}, \cite[Proposition~4.15]{CEVV}.
    Some examples of smoothable components other than the curvilinear component
    are given in~\cite[Appendix]{Michalek}.

    The problem below concerns local complete intersections. From the
    deformation theory point of view, the graded ones are very easy: deforming
    them means just deforming each equation separately. In the local case,
    things are more difficult, as deforming the equations can change the
    degree of the scheme. Things get even harder if we require that the scheme
    stays supported at a single point. The motivation for the problem below partially comes from
    motivic homotopy theory.
    \begin{problem}\label{prob:lci}
        In the punctual Hilbert scheme $\Hilb_d(\mathbb{A}^n,0)$ consider the
        open locus $\Hilb_d^{\lci}(\mathbb{A}^n, 0)$ of locally complete
        intersections supported at a point. Find an upper bound for
        its dimension.
    \end{problem}
    Locally complete intersections are smoothable, hence
    $\Hilb_d^{\lci}(\mathbb{A}^n, 0)$ is contained in
    $\Hilb_d^{\sm}(\mathbb{A}^n)$ and so it has dimension at most $nd-1$. However,
    local complete intersections are not necessarily limits of curvilinear
    schemes~\cite{ia_deformations_of_CI}\footnote{On page 602
        of~\cite{ia_deformations_of_CI} there is a bound $nd(1 - d^{-1/n})$
        attributed to Brian{\c{c}}on-Granger without specific reference. It would be a
        solution to Problem~\ref{prob:lci}. Unfortunately it seems wrong: the
        analogue of~\cite[Example~3]{ia_deformations_of_CI} for $C(4,4,4,4)$
        yields a family of too high dimension.}. Perhaps the newly introduced method of
        \BBname{} slicing~\cite{jelisiejew_keneshlou} could be of use for this
        problem.

    \subsection{Multigraded Hilbert schemes}

    Fix $n$ and $d$. Let $S = \kk[x_0, \ldots, x_n]$, note the indexing from
    zero.
    The multigraded Hilbert schemes are slightly different in flavour from the
    usual ones: they parameterize \emph{homogeneous} ideals $I$ in $S$ with
    $\dim(S/I) = 1$. This is still a ``zero-dimensional'' situation,
    when we intuitively look at these ideals as defining subschemes in $\mathbb{P}^n$ rather than
    in $\mathbb{A}^{n+1}$. We stress however, that the ideals above are
    usually \emph{not} saturated, hence there is important information lost
    when passing from $I$ to $V(I)\subseteq \mathbb{P}^n$.
    The multigraded Hilbert scheme is recently
    applied very successfully in theory of tensors, for investigating the lower
    bounds on border rank~\cite{Buczyska_Buczynski__border,
    Conner_Huang_Landsberg}.

    Let $H\colon \mathbb{N}\to \mathbb{N}$ be the Hilbert function of a tuple
    of $d$ general points in $\mathbb{P}^n$, so that $H(i) = \min(d,
    \binom{n+i}{i})$. Let $\Hilb^H$ be the multigraded Hilbert scheme as
    in~\cite{Haiman_Sturmfels__multigraded}, so its $\kk$-points correspond to
    homogeneous ideals $I\subseteq S$ with Hilbert function $H$. A bit more
    generally, let $H_r\colon \mathbb{N}^r\to \mathbb{N}$ be the Hilbert
    function of a tuple of $d$ general points in
    \[
        \underbrace{\mathbb{P}^{n} \times \ldots \times \mathbb{P}^n}_{r},
    \]
    then $\Hilb^{H_r}$ parameterizes homogeneous ideals in $I\subseteq S$
    where $S$ is an $\mathbb{N}^r$-graded ring in $r(n+1)$ variables
    corresponding to the total coordinate ring of $\mathbb{P}^{n} \times
    \ldots \times \mathbb{P}^n$.

    The multigraded Hilbert scheme $\Hilb^H$ behaves more pathologically than
    $\Hilbnd$. For $n=1$ it is smooth
    irreducible~\cite{Maclagan_Smith__multigraded}, but already for $n=2$ the
    multigraded Hilbert scheme is reducible and singular~\cite{Mandziuk},
    while $n=2$, na\"ively speaking, corresponds to projective plane and
    $\Hilb_d(\mathbb{P}^2)$ is smooth irreducible by Fogarty's result.

    \newcommand{\Satbar}{\overline{\operatorname{Sat}}}
    The saturated ideals in $\Hilb^H$ form an open subset~\cite{jelisiejew_Mandziuk}, whose closure is a
    union of irreducible components and is denoted by $\Satbar^H$. Inside,
    the saturated ideals of smooth subschemes in $\mathbb{P}^n$ form an open
    irreducible subset, whose
    closure is a distinguished component $\Slip^H$, the \emph{scheme of limits
    of ideals of points},
    see~\cite{Buczyska_Buczynski__border, Mandziuk, jelisiejew_Mandziuk}.
    The following problem is akin to smoothability and has immediate
    applications in the theory of border ranks of tensors:
    \begin{problem}
        For a given $[I]\in \Hilb^H$ determine whether it lies in~$\Satbar^H$
        or in $\Slip^H$.
    \end{problem}
    Determining this for monomial or even Borel-fixed ideals is already open.
    There are examples not lying in $\Satbar^H$, see~\cite{Mandziuk,
    jelisiejew_Mandziuk}.
    An extension from $H$ to $H_r$ would be also very interesting.
    See~\cite{jelisiejew_Mandziuk} for some results.

    By a result of Pardue~\cite[Theorem~34]{Pardue} the scheme $\Hilb^H$ is
    connected. The author is not aware of any extensions of this result to the
    case $r > 1$.
    \begin{problem}
        Is $\Hilb^{H_r}$ connected?  Does every component of $\Hilb^{H_r}$
        intersect the $\Slip^{H_r}$?
    \end{problem}

    If we remove the assumptions on $H_r$, that is, if we consider an arbitrary
    Hilbert function $H_r$ for an arbitrary $\mathbb{N}^r$-grading, then
    $\Hilb^{H_r}$ may be
    disconnected~\cite{Santos__Non_connected_toric_Hilbert_schemes}. The
    example is however complicated and the polynomial ring has $26$ variables.
    Hence we pose the following challenge.
    \begin{problem}
        Give a ``small'' example of a disconnected multigraded Hilbert scheme,
        preferably with $H_r$ eventually constant.
    \end{problem}

    \subsection{Smoothability} Smoothability of a given ideal is investigated
    classically. In this section we gather problems which were
    not posed before explicitly, however in our opinion they are very
    important.

    \begin{problem}\label{prob:defect}
        For a subscheme $\Gamma \subseteq \mathbb{A}^n$ consider the
        \emph{smoothability defect} $\delta(\Gamma)$ of $\Gamma$ as the minimal possible difference $\deg \Gamma' - \deg
        \Gamma$ where $\Gamma\subseteq \Gamma'$ and $\Gamma'$ is finite smoothable.
        What is the behavior of the function $\delta(-)$? For example, can it
        happen that $\delta(\Gamma) > \deg \Gamma $?
    \end{problem}
    As far as we know, nothing is known about the smoothability defect.
    To see that it is finite, take $\Gamma$ local and such that
    $|\Gamma| = \{0\}\subseteq \mathbb{A}^n = \Spec(\kk[x_1, \ldots ,x_n])$. Then $\Gamma$ is contained in
    $\Gamma' = V( (x_1, \ldots ,x_n)^e)$ for a large enough $e$. The scheme
    $\Gamma'$ is monomial, hence smoothable~\cite[Proposition~4.15]{CEVV}. For
    a nonlocal
    $\Gamma$ we repeat the above procedure for every point of the support of
    $\Gamma$.
    The above gives an exponential bound on the smoothability defect. However,
    in practice this defect seems to be very small, for example by a direct
    check it is equal to $1$ for
    the $(1,4,3)$ example of~\cite[Proposition~5.1]{CEVV} and the $(1,6,6,1)$ example of
    Emsalem-Iarrobino~\cite{emsalem_iarrobino_small_tangent_space}.

    \begin{problem}\label{prob:asymptotic}
        For a finite scheme $\Gamma$, let $\delta(\Gamma)$ be its smoothability defect.
        Is it true that for every $\epsilon > 0$ we have
        \[
            \lim_{n\to \infty} \frac{\delta(\Gamma^{\times n})}{(\deg
                \Gamma)^{n(1+\epsilon)}}
            = 0?
        \]
    \end{problem}

    Warning: an affirmative answer in Problem~\ref{prob:asymptotic}
    for every
    member of a sequence of $\Gamma := \Spec(\Abarn)$ as defined below
    (following~\cite{CHILO}) would
    imply the famous ``$\omega = 2$'' conjecture on the matrix
    multiplication exponent. This problem itself is a strengthened version
    of a special case of a conjecture by Strassen so \emph{disproving} it, even on an example
    disconnected from the context of the exponent of matrix multiplication, would be a breakthrough for
    complexity theory.
    See~\cite{Conner_Gesmundo_Landsberg_Vertura_Wang} for
    the discussion on Strassen's conjecture in the language of tensors.
    Providing an example of $\Gamma$ with nonzero limit for $\epsilon = 0$ is a bit
    weaker, but still very interesting problem.

    There is a more direct link between smoothability and the matrix
    multiplication exponent $\omega$.
    Consider a cubic form $\sMn$ in $n^2$ variables $(x_{ij})_{1\leq i,j\leq n}$
    given by $\tr([x_{ij}]^3)$, where $[x_{ij}]^3$ is the cube of a matrix
    $[x_{ij}]$. (The symbol $\sMn$ stands for
    \emph{symmetrized matrix multiplication tensor}~\cite{CHILO}.)
    The apolar algebra $(A_n,\mm)$ of $\sMn$
    is graded and has Hilbert function $(1,n^2, n^2, 1)$. Let $\Abarn =
    A_n/\mm^3$, then $\Abarn$ has Hilbert function $(1, n^2, n^2)$.
    Let $\delta_n := \delta(\Abarn)$ be the smoothability defect of $\Abarn$,
    as defined above. It follows from~\cite{CHILO} that $\delta_n$ grows
    like $n^{\omega}$ and that it is enough to look at surjections $B_n\onto
    \Abarn$ where $H_{B_n} = (1, *, *)$. (The details of this assertion will
    be given elsewhere.)
    Hence, to prove that $\omega = 2$
    it is enough to find for every $n$ a smoothable algebra of length around
    $n^2$ surjecting onto a given one! More modestly, it would be huge advance to
    find a sequence of $B_n$ which grows slower than, say, $n^{2.3}$.

    To formulate the next problem, let $\kk[[t]]$ be the ring of formal power series
    in one variable and $\kk((t))$ be the field of Laurent formal power
    series. For an element $f\in \kk((t))$ its \emph{valuation} is the
    smallest exponent $\nu(f)$ such that the coefficient in $f$ next to
    $t^{\nu(f)}$ is nonzero.

    Let $\kk((t))^{\times d} = \kk((t)) \times  \ldots \times
    \kk((t))$ be the algebra with usual coordinatewise operations.
    An \emph{abstract
    smoothing} of $A$ is a finite flat family $\kk[[t]]\to
    \mathcal{A}$ with generic fiber $\mathcal{A}_t$ (geometrically) reduced
    and special fiber $\mathcal{A}/t\mathcal{A}$ isomorphic to $A$.
    The scheme $\Gamma$ is smoothable in the usual sense if and only if an
    abstract smoothing
    exists~\cite[Theorem~1.1]{jabu_jelisiejew_smoothability}. A
    \emph{completely split abstract smoothing} of $A$ is an abstract smoothing
    $\mathcal{A}$ of $A$ with generic fiber $\mathcal{A}_{t}$ isomorphic to
    $\kk((t))^{\times d}$.

    For a finite flat $\kk[[t]]$-algebra $\mathcal{A}$,
    the \emph{approximation degree} of $\mathcal{A}$ is the maximum of
    valuations of nonzero elements appearing in the multiplication table of
    $\mathcal{A}$; this degree does not depend on the choice of a
    $\kk[[t]]$-basis.
    \begin{problem}\label{prob:approximation}
        Is it true that for every smoothable finite scheme $\Spec(A)$
        of degree $d$
        there exists a completely split abstract smoothing $\mathcal{A}$ of $A$ with
        approximation degree at most $d$?
        Less restrictively, does there exist a polynomial $P$ such that for
        every $A$ as above there exists an $\mathcal{A}$ with
        approximation degree bounded by $P(d)$? (for example, proving that
        the approximation degree is bounded by $d^{100}$ would be enough!)
    \end{problem}
    An abstract
    smoothing is rarely completely split. For example take $A = \kk[x]/(x^d)$.
    Then $\mathcal{A} = \kk[[t]][x]/(x^d - t)$ is an abstract smoothing of
    $A$, but $\mathcal{A}_t  \simeq \kk((t))[x]/(x^d-t)$ is a domain, so it is
    not isomorphic to $\kk((t))^{\times d}$. However, take
    \[
        \mathcal{A}' = \frac{\kk[[t]][x]}{(x(x-t)(x-2t) \ldots (x-(d-1)t))}.
    \]
    Assuming that $\kk$ has characteristic zero, by Chinese Remainder Theorem
    the generic fiber $\mathcal{A}'_t$ is indeed isomorphic to
    $\kk((t))^{\times d}$. The approximation degree obtained from
    $\mathcal{A}'$ is $d-1$. Observe that in the multiplication table of
    $\mathcal{A}$ the largest power of $t$ that appears is $t^1$, but this
    tells us nothing about the approximation degree since the smoothing $\mathcal{A}$ is not
    completely split. Completely split abstract smoothings exist for every
    smoothable $A$ and can be obtained by a base change from any
    one-parameter family that shows smoothability of $A$.

    Again, the motivation for the problem above is algebraic complexity
    theory~\cite{Ikenmeyer_talk, Lehmkuhl_Lickteig,
    Burgisser_complexity_factors}.
    In the references just given the problem is formulated differently. The
    author believes that the solution to Problem~\ref{prob:approximation}
    yields a polynomial upper bound on the approximation degree for tensors, a
    classical important direction in complexity theory~\cite{Lehmkuhl_Lickteig, Burgisser_complexity_factors,
    Burgisser_complexity_factors_erratum,
    Christandl_Gesmundo_Lysikov_Steffan}, but the precise argument on why it
    is so is nontrivial and well beyond the scope of the current work.

    \subsection{Moduli space of bilinear maps}\label{sec:bilinear}

    \newcommand{\Bilin}{\operatorname{Bilin}}%
    Let $S = \kk[x_1, \ldots ,x_n]$ be the coordinate ring of $\mathbb{A}^n$.
    Consider the \emph{Bilinear scheme} $\Bilin_{d_1, d_2, d_3}^{r_1, r_2}(\mathbb{A}^n)$ whose
    $\kk$-points are surjections of $S$-modules $S^{\oplus r_1}\onto
    M_1 = \frac{S^{\oplus r_1}}{K_1}$, $S^{\oplus
    r_2}\onto M_2 = \frac{S^{\oplus r_2}}{K_2}$ together with a surjection of $S$-modules
    $\pi\colon M_1\otimes_{S} M_2\onto M_3 = \frac{S^{\oplus
    r_1r_2}}{K_{3}}$, where $M_i$ is a finite $S$-module of
    degree $d_i$. This space can be realized as a scheme using Quot schemes.

    The following is a meta problem. Some interesting first results were obtained in~\cite{Obcowska_Bilinear}.
    \begin{problem}
        Investigate the geometry of this space.
    \end{problem}
    In practical terms, this means: find the components, investigate
    singularities, decide nonreducedness etc.
    Already
    the case $n = 2$, the case of two variables, is open.
    The motivation for this problem comes from the theory of tensors, but the
    Bilinear scheme naturally generalizes the Hilbert scheme and Quot scheme.
    To see this, we need to introduce a bit of notation.

    Observe that a
    finite algebra $A = S/I$ gives rise to a surjection $\pi_A\colon A \otimes_S A\to A$
    which is the multiplication in $A$. Every finite cyclic $S$-module is of the form
    $S^{\oplus 1}/I$ and a surjective $S$-linear map $\pi\colon S^{\oplus
    1}/I\otimes S^{\oplus 1}/J\to S^{\oplus 1}/K$ between modules with
    $\dim_{\kk} S/I = \dim_{\kk} S/J = \dim_{\kk} S/K$ exists only if $I = J =
    K$ and in that case it is equal to the  multiplication in $S/I$ followed
    by multiplying by a fixed invertible element of $S/I$; the latter
    operation is irrelevant when we pass to equivalence classes as in the Quot
    scheme.
    It follows that $\Bilin_{d_1, d_1, d_1}^{1,1}$ is
    isomorphic to
    $\Hilb_{d_1}(\mathbb{A}^n)$.
    More generally, the space $\Bilin_{d_1, d_1, d_1}^{r_1,r_2}$
    admits an open subset $\mathcal{U}$ where $M_1$ and
    $M_2$ are cyclic. On this open subset $M_3$ is cyclic as well (since it is an
    image of the cyclic module). This open subset is smoothly equivalent to
    the Hilbert scheme.

    For the Quot scheme, the situation is similar yet a bit more complicated,
    as we do not obtain exactly the Quot scheme: rather than modules, the
    obtained space parameterizes modules together with an algebra acting on
    them.
    Namely, a finite $S$-module $M = S^{\oplus
    r}/K$ yields a surjection $S/\Ann_S(M) \otimes_S M \to M$.
    The space $\Bilin_{d_1, d_2, d_2}^{r_1, r_2}$ admits an open subset
    $\mathcal{V}$
    where $M_1$ is cyclic. For a point of this subset, we have $M_1 \simeq
    S^{\oplus 1}/I$ and the homomorphism $\pi\colon M_1\otimes_S M_2\to M_3$
    simplifies to a surjective $S$-linear map $\frac{M_2}{IM_2}\to M_3$, which
    can only exist if $M_2  \simeq M_3$ and $IM_2 = 0$. Thus, from such a point
    we obtain a degree $d$ module $M := M_2  \simeq M_3$ and an algebra $A =
    S/I$ such that $IM = 0$, hence $M$ is an $A$-module. Let $\mathcal{Q}$ be the moduli space of pairs
    $M=S^{\oplus r_2}/K$ and $A = S/I$ such that $IM = 0$ and $\dim_{\kk} A =
    d_1$, $\dim_{\kk} M = d_2$. Then
    $\mathcal{V}$ admits a smooth map to
    $\mathcal{Q}$. Moreover, $\mathcal{Q}$ has a map to
    $\Quot_{d_2}^{r_2}(\mathbb{A}^n)$ which forgets about $A$. This map is usually not smooth. It may
    even be far from surjective: if $d_1$ is small, then the existence of $I$
    is very restrictive. If $d_1 = d_2$, then the map $\mathcal{Q}\to
    \Quot_{d_2}^{r_2}(\mathbb{A}^n)$ is birational on the principal component,
    hence $\mathcal{Q}$ can be thought of as a blowup of $\Quot_{d_2}^{r_2}(\mathbb{A}^n)$.

    In the language
    of~\cite{Jelisiejew_Landsberg_Pal__Concise_tensors_of_minimal_brank},
    roughly speaking,
    the Bilinear scheme corresponds to $111$-abundant tensors, the open subset
    $\mathcal{U}$ corresponds to $1_A$- and $1_B$-generic $111$-abundant
    tensors, while $\mathcal{V}$ corresponds to $1_A$-generic $111$-abundant
    tensors, see
    \cite[Proposition~5.3]{Jelisiejew_Landsberg_Pal__Concise_tensors_of_minimal_brank}.
    The passage between a tensor and bilinear map is straightforward: to a
    map $\pi\colon M_1\otimes_S M_2\to M_3$ one associates the tensor $\pi\in
    M_1^*\otimes_{\kk} M_2^* \otimes_{\kk} M_3$. Conversely, to a
    $111$-abundant tensor $T\in V_1\otimes_{\kk} V_2\otimes_{\kk}
    V_3$ with
    the $111$-algebra $A$, one associates the bilinear map $T\colon
    V_1^*\otimes_{\kk} V_2^*\to V_3$. The spaces $V_1^*$, $V_2^*$, $V_3$ are
    $A$-modules in a canonical
    way~\cite[\S5.1]{Jelisiejew_Landsberg_Pal__Concise_tensors_of_minimal_brank}
    and the map $T$ is $A$-bilinear, thus (after fixing some generators)
    yields a point of the Bilinear scheme\footnote{In~\cite{Jelisiejew_Landsberg_Pal__Concise_tensors_of_minimal_brank}
the spaces are named $A$, $B$, $C$ instead of $V_1$, $V_2$, $V_3$.}.

    There are also points of $\Bilin_{d_2, d_2,
    d_2}^{r_1, r_2}$ in the complement of $\mathcal{V}$.
    Examples for $d_1 = d_2 = 5$ are given
    in~\cite[Theorem~1.7]{Jelisiejew_Landsberg_Pal__Concise_tensors_of_minimal_brank}.
    Below, we translate the tensor called $T_{\OO_{58}}$ to the language
    above (see
    also~\cite[Example~4.6]{Jelisiejew_Landsberg_Pal__Concise_tensors_of_minimal_brank}).
    Here, we have $n = 4$, $d_1 = d_2 = d_3 = 5$ and $r_1 = r_2 =
    2$. We now construct the modules. Consider a polynomial ring in two variables $\kk[z_1, z_2]$ and the
    finite algebra $B = \frac{\kk[z_1, z_2]}{(z_1^4, z_2^4)}$.
    It has Hilbert function $H_B = (1, 2, 3, 4, 3, 2, 1)$. Let $B_j$ be the
    space of homogeneous elements of degree $j$ in $B$.
    Let $S = \kk[x_1,  \ldots , x_4]$ act on $B$ via the homomorphism
    $x_i\mapsto z_1^{i-1}z_2^{4-i}$ for $i=1,2,3,4$. Let $M$ be equal to the
    $\kk$-vector space $B_1 \oplus B_4$,
    so $M$ has a basis $z_1, z_2,
    z_1^3z_2, z_1^2z_2^2, z_1z_2^3$. Of course, $M\subseteq B$ is not a
    $B$-submodule. But $M\subseteq B$ \emph{is} an $S$-submodule.
    Let $N := B_{6-1} \oplus B_{6-4} = B_2\oplus B_5$. This is also an $S$-submodule of
    $B$. The multiplication in $B$ restricted to $M$ satisfies $M\cdot
    M\subseteq N$, so yields a map $\pi\colon
    M\otimes_S M\to N$. Choose $z_1, z_2$ as generators of the $S$-module $M$,
    so we have $M = S^{\oplus 2}/K$ and $\pi$ yields a point of
    $\Bilin_{5, 5, 5}^{2, 2}$. We stress that $B$ is not determined by the map
    $\pi$: we introduced it above as an ad-hoc device to phrase the example neatly.

\end{document}